\documentclass[a4paper,12pt]{article}

\usepackage{amsmath,amsfonts,amssymb,amsthm}
\theoremstyle{plain}
\newtheorem{thm}{Theorem}[section]

\newtheorem{lem}[thm]{Lemma}

\theoremstyle{definition}

\theoremstyle{remark}

\numberwithin{equation}{section}

\newcommand{\reftit}{\textit}    
\newcommand{\refis}{\textbf}     

\begin{document}

\title{Linear speed large deviations for percolation clusters}

\author{Yevgeniy Kovchegov and Scott Sheffield}

\maketitle \begin{abstract}

Let $C_n$ be the origin-containing cluster in subcritical percolation on the lattice $\frac{1}{n}
\mathbb Z^d$, viewed as a random variable in the space $\Omega$ of compact, connected,
origin-containing subsets of $\mathbb R^d$, endowed with the Hausdorff metric $\delta$.  When $d
\geq 2$, and $\Gamma$ is any open subset of $\Omega$, we prove that $$\lim_{n \rightarrow
\infty}\frac{1}{n} \log P(C_n \in \Gamma) = -\inf_{S \in \Gamma} \lambda(S)$$ where $\lambda(S)$ is
the one-dimensional Hausdorff measure of $S$ defined using the {\em correlation norm}: $$||u|| :=
\lim_{n \rightarrow \infty} - \frac{1}{n} \log P (u_n \in C_n )$$  where $u_n$ is $u$ rounded to
the nearest element of $\frac{1}{n}\mathbb Z^d$.  Given points $a^1, \ldots, a^k \in \mathbb R^d$,
there are finitely many correlation-norm Steiner trees spanning these points and the origin. We
show that if the $C_n$ are each conditioned to contain the points $a^1_n, \ldots, a^k_n$, then the
probability that $C_n$ fails to approximate one of these trees tends to zero exponentially in $n$.
\end{abstract}

\section{Introduction}

Let $C_n$ be the origin-containing cluster in subcritical Bernoulli bond-percolation with parameter
$p$ on the lattice $\frac{1}{n} \mathbb{Z}^d$; we view $C_n$ as a random variable in the space
$\Omega$ of compact, connected, origin-containing subsets of $\mathbb R^d$.  When the probability
measure involved is clear from context, we use $P(A)$ to denote the probability of an event $A$.
When $u \in \mathbb R^d$, let $u_n$ be the vector $u$ rounded to the nearest element in
$\frac{1}{n} Z^d$.  We define the ``correlation norm'' by  $$||u|| := \lim_{n \rightarrow \infty} -
\frac{1}{n} \log P (u_n \in C_n ).$$ This limit exists for all $u \in \mathbb R^d$ (with $||u|| \in (0,\infty)$
for $u \not = 0$) and $||\cdot||$ is a strictly convex norm (i.e., if $u$ and $v$ are not on
the same line through the origin, then $||u+v|| < ||u|| + ||v||$) that is real-analytic on the
Euclidean unit sphere $S^{d-1}$ \cite{ioffe}. Denote by $\lambda(S)$ the {\it one-dimensional Hausdorff measure} of the set $S$ defined with the above
norm; in particular, if $S \in \Omega$ is a finite union of rectifiable arcs in $\mathbb R^d$, then
$\lambda(S)$ is the sum of the correlation-norm lengths of those arcs.

Given a set $X \subset \mathbb R^d$, denote by $B_{\epsilon}(X)$ the set of all points of distance
less than $\epsilon$ from some point in $X$. Given sets $X, Y \in \Omega$, let $\delta(X,Y)$ be
their {\it Hausdorff distance}, i.e., $$\delta(X,Y) = \inf \{ \epsilon : X \subset B_{\epsilon}(Y),
Y \subset B_{\epsilon}(X) \}.$$

Many authors, including \cite{acc}, \cite{ccc}, \cite{dks}, \cite{ioffe}, and \cite{pisztora}, have
investigated the shapes of ``typical'' large finite clusters in supercritical percolation on
$\mathbb Z^d$ by proving surface order large deviation principles for clusters conditioned to
contain at least $m$ vertices.  They have shown that as $m$ gets large, the shapes of typical
clusters are approximately minimizers of surfaces tension integrals, called Wulff crystals.
Moreover, the surface tension integral is a rate function for a large deviation principle---with
{\it surface order} speed $m^{d-1/d}$---on cluster shapes.  These results are one way of precisely
answering the questions, ``What does the typical `large' cluster look like?  How unlikely are large
deviations from this typical shape?''

If instead of number of vertices we define ``large'' in terms of, say, diameter or volume of the
convex hull, then these questions can be answered for subcritical percolation using the following
linear speed large deviation principle: \begin{thm} \label{LDP} Let $d \geq 2$, $p < p_c$, and
$\Gamma \subset \Omega$ be Borel-measurable.  Then $$-\inf_{S \in \Gamma^o} \lambda(S) \leq
\liminf_{n \rightarrow \infty} \frac{1}{n}\log P(C_n \in \Gamma) \leq \limsup_{n \rightarrow
\infty}\frac{1}{n} \log P(C_n \in \Gamma) \leq -\inf_{S \in \overline{\Gamma}} \lambda(S)$$ where
$\Gamma^o$ and $\overline{\Gamma}$ are the interior and closure of $\Gamma$ with respect to the
Hausdorff topology. \end{thm}

In the language of \cite{dz}, this says that the random variables $C_n$ satisfy a large deviation
principle with respect to the Hausdorff metric topology on $\Omega$ and with speed $n$ and rate
function $I(S) = \lambda(S)$.  Note that since $\lambda: \Omega \rightarrow \mathbb R$ is
continuous, this implies that $$\lim \frac{1}{n}\log P(C_n \in \Gamma) =-\inf_{S \in \Gamma}
\lambda(S)$$ whenever $\Gamma$ is an open subset of $\Omega$.

\medbreak {\noindent\bf Acknowledgments.} We thank Amir Dembo and his probability discussion group
for helpful conversations and thank Yuval Peres for some suggestions on the presentation.  Also, Raphael
Cerf has informed us that Olivier Couronn\'{e}, working independently,
produced an alternate proof of Theorem \ref{LDP} and was nearly finished writing up the result at the time that
our paper was submitted and posted to the arXiv \cite{couronne}.

\section{Proof of large deviation principle} \label{proof}
\subsection{Exponential tightness and an equivalent formulation}
We now prove Theorem \ref{LDP}.  The sets $\{ S| \delta(S, \{0 \}) \leq \alpha \}$ are compact in
the Hausdorff metric topology, and $P(\delta(C_n, \{0\}) > \alpha)$ decays exponentially in $n$ and
$\alpha$. \cite{menshikov}  This implies that the laws of the $C_n$ are {\em exponentially tight}
(in the sense of \cite{dz}, Sec. 1.2). Given this exponential tightness, Theorem \ref{LDP} is
equivalent to the statement that the following bounds hold for $S \in \Omega$:

$$\lim_{\epsilon \rightarrow 0} \mathcal A(S,\epsilon) \leq \lambda(S)$$ $$\lim_{\epsilon
\rightarrow 0} \mathcal B(S,\epsilon) \geq \lambda(S)$$

where

 $$\mathcal A(S,\epsilon)=\limsup_{n \rightarrow \infty} \frac{-1}{n} \log P (\delta(S, C_n) < \epsilon )$$
 $$\mathcal B(S,\epsilon)=\liminf_{n \rightarrow \infty} \frac{-1}{n} \log P (\delta(S,
C_n) < \epsilon )$$ This equivalence is well-known in the large deviations literature (\cite{dz},
Lemma 1.2.18 and Theorem 4.1.11), and is also not hard to prove directly.  We now prove the first
of the two bounds above, which involves giving a lower bound on the probabilities $P (
\delta(S,C_n) < \epsilon )$.

\subsection{Lower bound on probabilities}
Fix $\epsilon$ and choose $S'$ to be a connected union of finitely many line segments of the form
$(a^i,b^i)$, for $1 \leq i \leq k$---intersecting one another only at endpoints---such that
$\delta(S,S') < \epsilon/2$ and at least one of the segments includes the origin as an endpoint. No
matter how small $\epsilon$ gets, we can always choose such an $S'$ of total length less than or
equal to $\lambda(S)$. Thus, it is enough to show that $$\liminf \frac{-1}{n} \log P ( \delta(S',
C_n) < \epsilon/2 ) \leq \lambda(S')$$ for sets $S'$ of this form.

Now, let $A^i_n$ (respectively, $A^i_{n,c}$ be the event that $a^i_n$ and $b^i_n$ are connected by {\it some} open path whose
Hausdorff distance from the line segment $(a^i, b^i)$ is at most $\epsilon/4$ (respectively $c/n)$.  For any
fixed $n$, $P(A^i_{n,c})$ tends to $P(a^i_n - b^i_n \in C_n)$ as $c$ tends to $\infty$.
Subadditivity arguments imply that $\liminf \frac{-1}{n} \log P(A^i_{n,c} )$ tends to $||a^i - b^i||$ as $c$ tends to
infinity.  It follows that $\liminf \frac{-1}{n} \log P(A^i_n ) \leq ||a^i - b^i||$.  
The FKG inequality then implies that $\liminf \frac{-1}{n} P ( \cup A^i_n ) \leq \lambda(S')$.

Now, we have to show that {\it given} $\cup A^i_n$, the probability of the event $C_n \not \subset
B_{\epsilon/2}(S')$ decays exponentially. Let $D_n$ be the event that there is a path from {\it
any} point $x$ outside of $B_{\epsilon/2}(S')$ to {\it any} point $y \in B_{\epsilon/4}(S')$. This
event is independent of $\cup A^i_n$. Since $D_n$ contains the event $C_n \not \subset
B_{\epsilon/2}(S')$, it is enough for us to show that $P(D_n)$ decays exponentially.  To see this,
we introduce and sketch a proof of the following lemma.  (See \cite{ioffe} for more delicate
asymptotics of $P ( u_n \in C_n )$.)

\begin{lem} \label{uniformbound} There exists a constant $\alpha$ such that $P ( u_n \in C_n ) \leq
\alpha e^{-n||u||}$ for all $n$ and $u$. \end{lem} \begin{proof} If $u=u_n$, then it is clear that
$P ( u_n \in C_n ) \leq e^{-n||u||}$.  (Simply use the FKG inequality to observe that for any
integer $m$, we have $P (u_{mn} \in C_{mn} ) \geq P (u_n \in C_n )^m$ and apply the standard
subadditivity argument to the log limits.)  If $u \not = u_n$, then it suffices to observe that
$e^{-n||u||}$ and $e^{-n||u_n||}$ differ by at most a constant factor. \end{proof}

The probability that any particular vertex of $B_{\epsilon}(S') \backslash B_{\epsilon/2}(S')$ is
connected to any particular vertex in $B_{\epsilon/4}(S)$ is bounded above by $\alpha \exp[ -n \inf
\{||u||: |u|=\epsilon/4\}]$, where $|u|$ is the Euclidean norm.  Since the number of pairs of
points of this type grows polynomially in $n$, the result follows.

\subsection{Upper bound on probabilities}
Fix $\gamma > 0$ and choose a finite set of points $a^1, a^2, \ldots a^k$ in $S$ such that every
collection $S'$ of line segments that contains the points $a^i$ has total length greater than
$\lambda(S)-\gamma$ (or greater than some large value $N$ if $\lambda(S)$ is infinite) and that for
some sufficiently small $\epsilon>0$, this remains true if each $a^i$ is replaced by some $c^i \in
B_{\epsilon}(a^i)$. (The reader may check that such a set of points and such an $\epsilon$ exist
for any $\gamma > 0$.) We know that $$\limsup -\frac{1}{n} \log P ( \delta(S, C_n)< \epsilon )$$ is
at least as large as $$\limsup -\frac{1}{n} \log P ( \mbox {some }c^i_n \in B_{\epsilon}(a^i)\mbox{
is contained in } C_n ).$$

We claim that the latter is at least $\lambda(S) - \gamma$. If $C_n$ does contain all of the
$c^i_n$, then it must contain a subgraph that is a tree with the $c^i_n$ as vertices. If we remove
all branches of this tree that do not contain a $c^i_n$, then a straightforward induction on $k$
shows that we are left with a tree $T$ in which at most $k-2$ vertices have more than two
neighbors. Denote by $b^i_n$ the vertices with this property. The path-connectedness-in-$T$
relation puts a tree structure on the set of $b^i_n$ and $c^i_n$. Each edge of this new tree $T'$
represents a pair of these points joined by a path, and all of these paths are disjoint.

Now, given a specific set of set of points $b^i_n$ and $c^i_n$ and $T'$, we have by the BK
inequality and Lemma \ref{uniformbound} that the probability that these disjoint paths are
contained in $C_n$ is at most $\alpha e^{-\lambda(T')n}$, where $\lambda(T')$ is the sum of the
correlation lengths of the edges of $T'$, and by assumption this value is at least $\lambda(S) -
\gamma$.  Since the number of possible choices for the $b^i_n$ and the $c^i_n$ grows polynomially,
and since $\gamma$ can be chosen arbitrarily small, this completes the proof.



\subsection{Steiner trees}
Given points $a^1,\ldots, a^k$, a (correlation norm) {\em Steiner tree spanning $\{a^i\}$ and the
origin} is an element $T$ of $\Omega$ for which $\lambda(T)$ is minimal among sets containing the
$\{a^i\}$. Existence of at least one Steiner tree follows from compactness arguments, and Steiner
trees are trees with at most $k-2$ vertices in addition to $a^1, \ldots, a^k$ \cite{gp}.  Although
the Steiner tree spanning a set of points is not always unique, strict convexity of the correlation
norm implies that the number of Steiner trees is always finite. See \cite{hrw} for a general
reference on Steiner trees. The proof of Theorem \ref{LDP} now yields the following:

\begin{thm} Let $d \geq 2$ and let $\Gamma \subset \Omega$ be Borel-measurable. If $C_n$ is the
origin-containing cluster in a subcritical percolation conditioned on $\{a^i_n\} \subset C_n$, then
$$-\inf_{S \in \Gamma^o} \lambda(S)- \lambda(T) \leq \liminf_{n \rightarrow \infty} \frac{1}{n}\log
P(C_n \in \Gamma) \leq \limsup_{n \rightarrow \infty}\frac{1}{n} \log P(C_n \in \Gamma) \leq
-\inf_{S \in \overline{\Gamma}} \lambda(S) - \lambda(T) $$ where $T$ is any Steiner tree spanning
$\{a^i\}$ and the origin. \end{thm}

In other words, these conditioned $C_n$ satisfy a large deviation principle with rate function
given by $I(S) = \lambda(S) - \lambda(T)$.  In particular, if $T_j$, for $1 \leq j \leq m$, are the
Steiner trees spanning $\{a^i\}$ and the origin, and $B_{\epsilon}(T_j) = \{S : \delta(S, T_j) <
\epsilon \}$, then we have

$$\lim_{n \rightarrow \infty}\frac{1}{n} \log P(C_n \not \in \cup B_{\epsilon}(T_j) = -\inf_{S \not
\in \cup B_{\epsilon}(T_j)} \lambda(S) - \lambda(T).$$  That is, the probability that $C_n$ fails to
approximate one of these Steiner trees tends to zero exponentially.

\bibliographystyle{amsplain}

\end{document}